\documentclass[12pt]{amsart}
\usepackage{amsmath}

\newcommand{\nc}{\newcommand}
\nc{\cal}{\mathcal} 
\nc{\la}{\langle}
\nc{\ra}{\rangle}
 \nc{\CA}{\cal A}
 \nc{\CBB}{\cal B}
 \nc{\CC}{\cal C}
\nc{\CDD}{\cal D}
\nc{\CE}{\cal E}
\nc{\CF}{\cal F}
\nc{\CG}{\cal G}
\nc{\CH}{\cal H}
\nc{\CI}{\cal I}
\nc{\CJ}{\cal J}
\nc{\CK}{\cal K}
\nc{\CL}{\cal L}
\nc{\CM}{\cal M}
\nc{\CN}{\cal N}
\nc{\CO}{\cal O}
\nc{\CP}{\cal P}
\nc{\CQ}{\cal Q}
\nc{\CR}{\cal R}
\nc{\CS}{\cal S}
\nc{\CT}{\cal T}
\nc{\CU}{\cal U}
\nc{\CV}{\cal V}
\nc{\CW}{\cal W}
\nc{\CZ}{\cal Z}

\nc{\fa}{\mathfrak a}
\nc{\fg}{\mathfrak g}
\nc{\fk}{\mathfrak k}
\nc{\fh}{\mathfrak h}
\nc{\fm}{\mathfrak m}
\nc{\fn}{\mathfrak n}
\nc{\fA}{\mathfrak A}
\nc{\fC}{\mathfrak C}
\nc{\fI}{\mathfrak I}
\nc{\fL}{\mathfrak L}
\nc{\fS}{\mathfrak S}

\nc{\nen}{\newenvironment}
\nc{\ol}{\overline}
\nc{\ul}{\underline}
\nc{\lra}{\longrightarrow}
\nc{\lla}{\longleftarrow}
\nc{\Lra}{\Longrightarrow}
\nc{\Lla}{\Longleftarrow}
\nc{\Llra}{\Longleftrightarrow}
\nc{\hra}{\hookrightarrow}
\nc{\iso}{\overset{\sim}{\lra}}
\nc{\Hom}{\mathrm{Hom}}
\nc{\Mor}{\mathrm{Mor}}

\nc{\notebox}[1]{\noindent\fbox{\parbox{12.5cm}{\sf #1}}\\[8pt]}
\nc{\Thm}[1]{Theorem~\ref{#1}}
\nc{\Prop}[1]{Proposition~\ref{#1}}
\nc{\Lem}[1]{Lemma~\ref{#1}}
\nc{\Cor}[1]{Corollary~\ref{#1}}
\nc{\Conj}[1]{Conjecture~\ref{#1}}
\nc{\Claim}[1]{Claim~\ref{#1}}
\nc{\Defn}[1]{ Definition~\ref{#1}}
\nc{\Exa}[1]{Example~\ref{#1}}
\nc{\Rem}[1]{Remark~\ref{#1}}
\nc{\Note}[1]{Note~\ref{#1}}

\nc{\marg}{\marginpar}

\nen{thm}[1]{\label{#1}{\bf Theorem.\ } \em}{}
\nen{prop}[1]{\label{#1}{\bf Proposition.\ } \em}{}
\nen{lem}[1]{\label{#1}{\bf Lemma.\ } \em}{}
\nen{klem}[1]{\label{#1}{\bf Key Lemma.\ } \em}{}
\nen{cor}[1]{\label{#1}{\bf Corollary.\ } \em}{}
\nen{conj}[1]{\label{#1}{\bf Conjecture.\ } \em}{}
\nen{claim}[1]{\label{#1}{\bf Claim.\ } \em}{}

\nen{defn}[1]{\label{#1}{\bf Definition.\ } }{}
\nen{exa}[1]{\label{#1}{\bf Example.\ } }{}

\nen{rem}[1]{\label{#1}{\em Remark.\ } }{}
\nen{exer}[1]{\label{#1}{\em Exercise.\ } }{}
\nen{pf}{\begin{proof}}{\end{proof}}

 \nc{\br}{\mathbb{R}}
 \nc{\bz}{\mathbb{Z}}
 \nc{\bc}{\mathbb{C}}
 \nc{\bn}{\mathbb{N}}
 \nc{\geg}{\mathfrak g}
 \nc{\G}{\Gamma}
 \nc{\sm}{\setminus}
 \nc{\sub}{\subset}
 \nc{\lm}{\lambda}
 \nc{\al}{\alpha}
  \nc{\bt}{\beta}
 \nc{\om}{\omega}
  \nc{\dl}{\delta}
  \nc{\g}{\gamma}
 \nc{\Dl}{\Delta}
 \nc{\Om}{\Omega}
 \nc{\s}{\sigma}
 \nc{\ro}{\rho}
  \nc{\te}{\theta}
 \nc{\SLR}{SL_2(\br)}
 \nc{\GLR}{GL_2(\br)}
 \nc{\PGLR}{PGL_2(\br)}
 \nc{\PSLR}{PSL_2(\br)}
 \nc{\SLC}{SL(2,\bc)}
 \nc{\uH}{\mathbb H}
 \nc{\fD}{\mathfrak D}
 \nc{\fE}{\mathfrak E}
 \nc{\haf}{\frac{1}{2}}
 \nc{\qtr}{\frac{1}{4}}
 \nc{\8}{\infty}
 \nc{\7}{{-\infty}}
 \nc{\inv}{^{-1}}
 \nc{\eps}{\varepsilon}

\begin{document}
\title[microlocal lifts]
{Microlocal lifts of eigenfunctions on hyperbolic surfaces and
trilinear invariant functionals}


\author{Andre Reznikov}
\address{Bar-Ilan University, Ramat Gan, Israel}
\email{reznikov@math.biu.ac.il}

\begin{abstract} In \cite{Z1} S. Zelditch introduced an
equivariant version of a pseudo-differential calculus on a
hyperbolic Riemann surface. We recast his construction in terms of
trilinear invariant functionals on irreducible unitary
representations of $\PGLR$. This allows us to use certain
properties of these functionals in the study of the action of
pseudo-differential operators on eigenfunctions of the Laplacian
on hyperbolic Riemann surfaces.
\end{abstract}
\maketitle

\section{Introduction}
\label{intro}
\subsection{Motivation}\label{motivation}

Let $Y$ be a compact Riemann surface  with a Riemannian metric of
constant curvature $-1$ and the associated volume element $dv$.
The corresponding  Laplace-Beltrami operator $\Dl$ is non-negative
and has purely discrete spectrum on the space $L^2(Y,dv)$ of
functions on $Y$. We will denote by $0=\mu_0< \mu_1 \leq \mu_2
\leq ...$ the eigenvalues of $\Dl$ and by $\phi_i=\phi_{\mu_i}$
the corresponding eigenfunctions (normalized to have $L^2$ norm
one). In the theory of automorphic functions the functions
$\phi_{\mu_i}$ are called non-holomorphic forms, Maass forms
(after H. Maass, \cite{M}) or simply automorphic functions.

The study of Maass forms is important in analysis and in other
areas. It has been understood since the seminal works of A.
Selberg \cite{Se} and I. Gel'fand, S.Fomin \cite{GF} that
representation theory plays an important role in this study.
Central for this role is the correspondence between eigenfunctions
of Laplacian on $Y$ and unitary irreducible representations of the
group $\PGLR$ (or what is more customary of $\PSLR$). This
correspondence allows one, quite often, to obtain results that are
more refined than similar results for the general case of a
Riemannian metric of variable curvature.

A framework where the correspondence between eigenfunctions and
representations plays a decisive role is the equivariant
pseudo-differential calculus constructed by S. Zelditch in
\cite{Z1}-\cite{Z3}. His motivation was to give a proof of the
celebrated quantum ergodicity theorem of A. Shnirelman for
hyperbolic surfaces (see Shnirelman \cite{Sh}, Y. Colin de
Verdi\`{e}re {\cite{CdV}, Zelditch \cite{Z2}). The main ingredient
of the proof of quantum ergodicity is a construction for each
eigenfunction $\phi_i$ on $Y$ of an associated probability measure
$dm_i$ on the spherical bundle $S^*(Y)$ of the co-tangent bundle
of $Y$. The idea to associate such measures to eigenfunctions was
a deep insight of Shnirelman. The measures $dm_i$ are called
microlocal lifts or micro-localizations of the corresponding
eigenfunctions $\phi_i$. The main property of these measures is
that they satisfy the Egorov-type theorem, that is, the measures
$dm_i$ are asymptotically invariant under the geodesic flow as
$\mu_i\to\8$. The measures $dm_i$ are constructed in two steps.
First, one constructs the so-called Wigner distribution $dU_i$  or
the auto-correlation distribution corresponding to $\phi_i$.
Namely, the distribution $dU_i\in\CDD(S^*(Y))$ such that for any
pseudo-differential operator (PDO) $A$ of order $0$ with the
symbol $a\in C^\8(S^*(Y))$ the relation
$<A\phi_i,\phi_i>=\int_{S^*(Y)}a\ dU_i$ holds. The distribution
$dU_i$ depends on a choice of the pseudo-differential calculus on
$Y$. Next, one modifies $dU_i$ (which are not non-negative) in
order to get a probability measure $dm_i$ asymptotically close to
$dU_i$ as $\mu_i\to\8$. Such a modification is not unique. The
Quantum Ergodicity Theorem claims that the measures $dm_i$
converge to the standard Liouville measure on $S^*(Y)$, at least
along a sequence of full density.

\subsection{Results}\label{results}
In this paper we discuss a relation of the measures $dm_i$ from
Zelditch's version of the equivariant pseudo-differential calculus
on $Y$ to representation theory of the group $\PGLR$. Namely, we
will show how the asymptotic invariance under the geodesic flow of
these measures follows from the uniqueness of invariant trilinear
functionals on three irreducible unitary representations of
$\PGLR$. This is based on the following theorem which is the main
underlying observation of the paper and comes from the uniqueness
of trilinear functionals and Zelditch's description of
pseudo-differential operators. To state it we recall first some
basic facts about $S^*(Y)$ and differential operators on this
space (see the standard excellent source \cite{GGPS}).

It is well-known (and fundamental for our approach) that there is
a transitive action of the group $G=\PGLR$ on $S^*(Y)$(usually one
considers the action of the group $\PSLR$, but for some technical
reasons explained in Section \ref{rem-uni}, we prefer to work with
$\PGLR$). Let $\mathcal{C}$ be the Casimir operator acting on the
space $C^\8(S^*(Y))$ of smooth functions on $S^*(Y)$ (this is the
unique up to a constant second order hyperbolic $G$-invariant
differential operator). The set of eigenvalues of $\mathcal{C}$
coincides with the set of eigenvalues of $\Dl$ on $Y$ (although
eigenspaces of $\mathcal{C}$ are infinite dimensional). Let
$V_\mu$ be a $\mu$-eigenspace of $\mathcal{C}$ which is {\it
irreducible} under the $G$-action. For an eigenvalue $\mu$ the
$\mu$-eigenspace splits into a direct sum of finitely many
irreducible ones and their span over all $\mu_i$ is dense in
$C^\8(S^*(Y))$ (see \ref{decompos}). It turns out that the space
$V_\mu$ is a unitarizable irreducible representation of $G$. The
representation $V_\mu$ is called an automorphic representation.
All unitarizable irreducible representation of $G$ are classified
(see \cite{G5} and Section \ref{reps} below). It turns out that
any irreducible unitarizable representation has a dense subspace
(called a space of smooth vectors) which could be realized as a
quotient of the space $C^\8(S^1)$ by a finite-dimensional
subspace. Hence for any space $V_\mu$ there exists a map
$\nu_\mu:C^\8(S^1)\to V_\mu\subset C^\8(S^*(Y))$. We assume for
simplicity that $\nu_\mu$ has no kernel and that $\mu\geq \qtr$.
It turns out that in this case the map $\nu_\mu$ gives rise to an
isometry $L^2(S^1)\to L^2(S^*(Y)$) and we denote by $<,>_{S^1}$,
$<,>_Y$ the corresponding scalar products and the corresponding
pairing between distributions and functions. The representation
$V_\mu$ with the above property is called a class one
representation of $G$ of principal series; for $V_\mu$ which is
not of class one the map $\nu_\mu$ has a finite dimensional
non-zero kernel. Such a representation $V_\mu$ is called a
discrete series representation. We deal with these also (see
\ref{main-proof}). Hence, any eigenfunction (in $V_\mu$) of
$\mathcal{C}$ is of the form $\nu_\mu(v)$ for an appropriate
function $v\in C^\8(S^1)$.

Let $Op$ be the pseudo-differential calculus of Zelditch (we
recall briefly the construction of $Op$ in Section
\ref{pdo-sect}). In particular, this calculus assigns to a symbol
$a(x,\lm)\in C^\8(S^*(Y)\times\br^+)$ an operator $Op(a)$ acting
on $C^\8(Y)$.

We will be interested in pseudo-differential operators of order
$0$ and moreover in those with the symbol independent of $\lm$
(see \ref{pdo-subsect}). This means that we consider the
correspondence between symbols $a\in C^\8(S^*(Y))$ and operators
$Op(a)$ acting on $C^\8(Y)$. For Maass forms Zelditch found a
description of the action of such pseudo-differential operators in
terms of the Helgason transform (which is a non-Euclidian analog
of the Fourier transform). We rephrase his description in terms of
representation theory as follows. Let $\mu$ be an eigenvalue of
$\Dl$ and $E_\mu\subset C^\8(Y)$ the corresponding eigenspace and
let $W_\mu$ be the corresponding eigenspace of $\mathcal{C}$ (we
have $W_\mu\simeq E_\mu\otimes V_\mu$). It turns out that one can
construct a map $\mathcal{M}:E_\mu\to W^*_\mu\subset
\mathcal{D}(S^*(Y))$ (called microlocalization) from the space of
Maass forms to the space of distributions on $S^*(Y)$ such that
for any symbol $a$ and any Maass form $\phi\in E_\mu$ we have
$Op(a)\phi(z)=\int_{S^1}a(zt)\mathcal{M}(\phi)(zt)dt$, where $z\in
Y$ and the integration is along the fiber of $S^*(Y)\to Y$ (i.e.
$Op(a)\phi$ is the push-forward of the distribution
$a\mathcal{M}(\phi)$ to $Y$. It turns out that the result is a
smooth function on $Y$ and hence the integration is well-defined).
Using this interpretation, we can express the action of a
pseudo-differential operator on an eigenfunction as the
multiplication of the corresponding distribution by the (smooth)
symbol of the operator. This allows us to relate
pseudo-differential operators to multiplication of automorphic
functions and then to trilinear invariant functionals on
representations.

We state now our main theorem (see \ref{main-proof})

\begin{thm}{main-thm} Let $V_\mu\subset C^\8(S^*(Y))$
be an irreducible eigenspace and $\nu_\mu:C^\8(S^1)\to V_\mu$ the
corresponding map and let $\mu_1$, $\mu_2$ be eigenvalues of the
Laplacian $\Dl$. There exists an explicit distribution
$l_{\al,\beta,\gamma}\in \mathcal{D}(S^1)$ on $S^1$ depending on
three complex parameters $\alpha,\beta,\gamma\in\bc$ such that for
any symbol $a$ of the form $a=\nu_\mu(v_a)\in C^\8(S^*(Y))$ with
$v_a\in C^\8(S^1)$ and $\phi_1$, $\phi_2$ eigenfunctions  of the
Laplacian $\Dl$, $\Dl\phi_i=\mu_i\phi_i$, there exists a constant
$a_{\mu,\mu_1,\mu_2}\in\bc$ satisfying the relation
\begin{align}\label{Op-tripl}
<Op(a)\phi_1,\phi_2>_Y=a_{\mu,\mu_1,\mu_2}\cdot
<l_{\mu,\mu_1,\mu_2},v_a>_{S^1} \ .
\end{align}
\end{thm}

Hence for the special kind of symbols, which we call irreducible
symbols (i.e. those belonging to one of the irreducible
representations $V_\mu$) we are able to analyze the action of the
corresponding pseudo-differential operator on eigenfunctions by
means of representation theory. We note that the space spanned by
such symbols is dense in the space of all smooth symbols.

We want to stress that for us the most important conclusion of the
theorem above is the claim that the distribution
$l_{\mu,\mu_1,\mu_2}$ has an {\it explicit} kernel which depends
only on eigenvalues and not on the choice of eigenfunctions
$\phi_i$ nor on the choice of the {\it symbol} $a$. We will use
this heavily throughout the paper. We will see that one can choose
the kernel of $l_{\al,\beta,\gamma}$ to be given by a function on
$S^1$ which is similar to the function
$|\sin(\theta)|^{\frac{-1-\lm}{2}}$ with $\lm$ pure imaginary (see
\eqref{integ-reduc1}).

The coefficients $a_{\mu,\mu_i,\mu_j}$ depend on the choice of
$\phi_i$ and $\phi_j$ and encode an important information about
corresponding eigenfunctions (for arithmetic surfaces and a
special basis of eigenfunctions, the Hecke-Maass basis, these
coefficients are connected to special values of certain
$L$-functions, see \cite{Sa2}). We will discuss bounds on these
coefficients as functions of eigenvalues $\mu_i$ and make some far
reaching conjectures about their size (see \ref{coef-prop}).

In Section \ref{3prod-def} we show how to express the setting of
pseudo-differential operators in terms of trilinear invariant
functionals on irreducible representations of $G$. The main
technical fact about trilinear invariant functionals we use in
this paper, beside their uniqueness, is that such functionals
could be described in terms of an explicit kernel. We study this
kernel from the point of view of oscillatory integrals.  Once we
relate the distribution $l_{\mu,\mu_i,\mu_j}$ to a trilinear
functional, we are able to give an explanation for the asymptotic
invariance of the microlocal measures in terms of the geometry of
the phase of this kernel. We also explain why some probability
measures suggested by the construction of S. Wolpert (\cite{Wo})
are asymptotic corrections to distributions $dU_i$. This gives the
positivity result necessary in the Shnirelman's argument.

We note that invariant trilinear functionals play an important
role in \cite{Z3}, albeit implicitly. Essentially, different
iterative formulas in \cite{Z3} (which were developed in order to
prove the asymptotic invariance in the first place) follow from
the uniqueness of invariant trilinear functionals (we note that
these formulas served as a starting point for the recent approach
of E. Lindenstrauss to the quantum unique ergodicity, see
\cite{Li}). The approach taken in \cite{Z3} is based on
differential relations coming from the action of the Lie algebra
$sl_2(\br)$ while our approach is based on properties of integral
operators involved and hence, in principle, is more flexible.
While the uniqueness of invariant trilinear functionals is widely
known to specialists in automorphic functions (where it plays an
important role in the theory of $L$-functions) it is rarely used
by analysts and deserves a wider recognition (see the recent book
\cite{U} however).

The paper is organized as follows. We begin with a brief review of
Zelditch's construction of the equivariant pseudo-differential
calculus (see \cite{Z1} for more detail) and the well-known
relation between eigenfunctions on $Y$ and representation theory
of $\PGLR$ (due to Gel'fand and Fomin, see \cite{GF},
\cite{GGPS}). We also express Zelditch's distributions $dU_i$ in
terms of special vectors (distributions) in the corresponding
automorphic representations. We then introduce our main tool of
invariant trilinear functionals and recast pseudo-differential
calculus in these terms. We next describe invariant trilinear
functionals explicitly in terms of their kernels. Central for this
is the alluded above uniqueness of invariant trilinear
functionals. It turns out that one can choose such kernels to be
given by simple homogeneous functions on copies of $\br^2\sm 0$.
To see this we use the standard model of the irreducible unitary
representations of $\PGLR$ realized in homogeneous functions on
$\br^2\sm 0$. We use this explicit description of trilinear
functionals in order to deduce the asymptotic invariance of
microlocal lifts of eigenfunctions (Theorem \ref{invar-thm}) and
to construct asymptotic probability measures (Theorem
\ref{asym-dm}). Both results follow from the explicit form of the
kernel of trilinear functional and the stationary phase method. We
also give a quantitative bound on the non-invariant part in terms
of an appropriate Sobolev norm. On the basis of our analysis we
show that one might expect that matrix coefficients
$<A\phi_i,\phi_i>$ are invariant under the geodesic flow up to a
higher order (by the factor $\mu_i^{-\frac{1}{4}}$) than the
Egorov's theorem predicts (this was also noticed by Zelditch). We
also show that for a fixed pseudo-differential operator $A$ the
spectral density of $A\phi_i$ is (essentially) supported in a
short interval near $\mu_i$ (Theorem \ref{prop}) and formulate a
conjecture concerning the size of coefficients $<A\phi_i,\phi_j>$
on this interval.


{\bf Acknowledgments.} This paper is a part of a joint project
with J. Bernstein whom I would like to thank for numerous fruitful
discussions. I would like to thank L. Polterovich for helpful
remarks which led to an improvement of the exposition.

The research was partially supported by BSF grant, Minerva
Foundation and by the Excellency Center ``Group Theoretic Methods
in the Study of Algebraic Varieties'' of the Israel Science
Foundation, the Emmy Noether Institute for Mathematics (the Center
of Minerva Foundation of Germany).

\section{Equivariant pseudo-differential operators}\label{pdo-sect}
We describe the construction of Zelditch \cite{Z1} of the
equivariant pseudo-differential calculus on a hyperbolic surface.
It is based on Helgason's representation theorem for
eigenfunctions on the unit disk $D$.

\subsection{Geometric setting}\label{geo-setting} We begin with
some well-known definitions (\cite{He},\cite{Z1}). Let $D$ be the
Poincar\'{e} open unit disk with the hyperbolic metric
$ds^2=(dx^2+dy^2)/(1-r^2)^2$ and the hyperbolic volume element
$dvol_H=dxdy/(1-r^2)^2$, where $r^2=x^2+y^2$. We denote by $B$ the
boundary circle of $D$ (on infinity) . Given a pair $(z,b)\in
D\times B$ let $\xi(z,b)$ be the unique horocycle through $z\in D$
with forward end point $b\in B$. The non-Euclidian (signed)
distance from the origin $0$ to $\xi(z,b)$ will be denoted
$<z,b>$. It is well known that functions
$e^{(\frac{\lm-1}{2})<z,b>}$ are eigenfunctions of the hyperbolic
Laplacian with the eigenvalue $\mu=\frac{1-\lm^2}{4}$ (here we
slightly changed normalization from the one adopted in \cite{Z1}).
The group $PSU(1,1)\simeq P\SLR$ acts by the standard fractional
linear transformations on $D$ and coincides with group of
isometries of $D$. We will use the identification $D\times B\simeq
P\SLR$ via the equivariant map sending a pair $(z,b)$ to the
unique element $g_{z,b}\in P\SLR$ such that $g_{z,b}\cdot 0=z$ and
$g_{z,b}\cdot 1=b$. One can view this as a well-known
identification $P\SLR\simeq S^*(D)\simeq S(D)$ with the
(co-)spherical bundle on $D$. The action of $g\cdot (z,b)\to
(gz,gb)$ coincides then with the left action of $P\SLR$ on itself.

We choose  $\G\subset P\SLR$ a (co-compact) discrete subgroup such
that the Riemann surface $Y=\G\simeq D$.

\subsection{Helgason's representation}\label{He-reps-sect} In
\cite{He} Helgason proved the following

\begin{thm}{He-thm} Let $\phi\in C^\8(D)$ of at most polynomial
growth (in the hyperbolic distance from the origin) near the
boundary $B$ and satisfying $\Dl\phi=\frac{1-\lm^2}{4}\phi$. Then
there exists a distribution on the boundary $T\in \CDD(B)$ such
that
\begin{align}\label{He-int}
\phi(z)=\int_B e^{(\frac{1+\lm}{2})<z,b>}dT(b)\ .
\end{align}
\end{thm}
We denote the correspondence defined by \eqref{He-int} by
$\CP_\lm:\CDD(B)\to C^\8(D)$ and refer to it as the Helgason map
(it is also called the non-Euclidian Poisson map). An important
point is that Helgason's representation is equivariant with
respect to the standard action of $P\SLR$ on $D$ and the following
twisted action on $B$. Namely, let $\pi_\lm$ be the representation
of $\SLR$ on the space of functions (or distributions) defined by
\begin{align}\label{lm-act-B}
\pi_\lm(g)f(b)=f(g\inv\cdot b)|g'(b)|^{\frac{\lm-1}{2}}\ .
\end{align}
This defines a representation of $P\SLR$ (which is unitary and
irreducible in the space $L^2(B)$ for $\lm\in i\br$). We have then
$\CP_\lm(\pi_\lm(g)T)(z)=\CP_\lm(T)(g\inv z)$ for any $g\in
P\SLR$. In particular, an eigenfunction $\phi$ is $\G$-invariant
if and only if the distribution $T$ is $(\pi_\lm,\G)$-invariant
(we will see latter that this is exactly the Frobenius reciprocity
from the theory of automorphic functions; see \ref{frob-sect}).

We note that there is an inverse to $\CP$ map given by (properly
defined) boundary values of eigenfunctions (see
\cite{He},\cite{Le}).

Similar to (\ref{He-int}) one have the following Helgason
non-Euclidian Fourier transform $\CF:C_0^\8(D)\to C^\8(\br^+\times
B)$ for a general function $f\in C_0^\8(D)$:
\begin{align}\label{Four-transf}
\CF(f)(\lm,b)=\hat
f(\lm,b)=\int_De^{(\frac{1-\lm}{2})<z,b>}f(z)dvol(z)\
\end{align}
and the inverse transform
\begin{align}\label{Four-inv-transf}
f(z)=\hat f(\lm,b)=\int_{\br^+\times
B}e^{(\frac{1+\lm}{2})<z,b>}\hat f(\lm,b)\lm\tanh(\pi\lm/2)d\lm
db\ .
\end{align}
The non-Euclidian Fourier transform $\CF$ is an isometry between
spaces $L^2(D,dvol_H)$ and $L^2(\br^+\times
B,(1/2\pi)\tanh(\pi\lm/2)d\lm db)$.

\subsection{Pseudo-differential operators}\label{pdo-subsect}
Based on the representation (\ref{Four-inv-transf}) Zelditch
introduced in \cite{Z1} the following form of $\SLR$-equivariant
pseudo-differential calculus.

Given any operator $A:C^\8(D)\to C^\8(D)$ one defines its complete
symbol $a(z,\lm,b)\in C^\8(D\times\br^+\times B)$ by
\begin{align}\label{symbol-def}
Ae^{ (\frac{1+\lm}{2})<z,b>}=a(z,\lm,b)e^{(\frac{1+\lm}{2})<z,b>}\
.
\end{align}
By the inversion formula (\ref{Four-inv-transf}), we have the
following representation
\begin{align}\label{operA-def}
Af(z)=1/2\pi\int_{\br^+\times
B}e^{(\frac{1+\lm}{2})<z,b>}a(z,\lm,b)\hat
f(\lm,b)\lm\tanh(\pi\lm/2)d\lm db\ .
\end{align}
It is assumed that the symbol of $A$ has the standard asymptotic
(in the symbol topology) expansion $a\sim
\sum_0^\8\lm^{-j}a_{-j}(z,b)$ as $|\lm|\to\8$.
 We will be interested in pseudo-differential operators of order
$0$ and hence will assume that the symbol is independent of $\lm$.
For such a symbol $a(z,b)\in C^\8(D\times B)$ we will denote
 $Op(a)$ the pseudo-differential operator defined by
 (\ref{operA-def}).

The correspondence between operators $A$ and their symbols
$a(z,\lm,b)$ is equivariant. Namely, the symbol of $gA$ is given
by $a(gz,\lm,gb)$. We will be interested in $\G$-invariant version
of pseudo-differential operators, i.e. those which commute with
the action of $\G$. Such symbols naturally gives rise to the
pseudo-differential operators on the Riemann surface $Y$.  Let
$\phi\in C^\8(\G\sm D)\simeq C^\8(Y)$ be an eigenfunction of $\Dl$
with the eigenvalue $\mu=\frac{1-\lm^2}{4}$ and $T\in\CDD(B)$ be
the boundary distribution assigned to $\phi$ via Helgason's
representation (\ref{He-int}). Zelditch then showed in \cite{Z1}
that for any $\G$-invariant symbol $a(z,b)\in C^\8(D\times B)^\G$
we have as above
\begin{align}\label{Op-inv-def}
Op(a)\phi(z)=\int_B a(z,b)e^{(\frac{1+\lm}{2})<z,b>} dT(b)\ .
\end{align}
This formula will serve us as a starting point for an
interpretation of $Op(a)$ in terms of representation theory and
particularly in terms of trilinear invariant functionals.

\section{Representation theory and eigenfunctions} \label{reps}
We recall the standard connection between eigenfunctions and
representation theory (see \cite{GGPS}).

\subsection{Automorphic representations} \label{ureps}

Let us describe the geometric construction which allows one to
pass from analysis on a Riemann surface to representation theory.

One stars with the Poincar\'{e} unit disk $D$ as above (or
equivalently $\uH$ the upper half plane with the hyperbolic metric
of constant curvature $-1$; the use of $\uH$ is more customary in
the theory of automorphic functions). The group $\SLR\simeq
SU(1,1)$ acts on $D$ (or $\uH$) by fractional linear
transformations. This action allows to identify the group $\PSLR$
with the group of all orientation preserving motions of $D$. For
reasons explained bellow we would like to work with the group $G$
of all motions of $D$; this group is isomorphic to $\PGLR$. Hence
throughout the paper we denote $G=\PGLR$.

Let us fix a discrete co-compact subgroup $\G \subset G$ and set
$Y=\G \sm D$. We consider the Laplace operator on the Riemann
surface $Y$ and denote by $\mu_i$ its eigenvalues and by $\phi_i$
the corresponding normalized eigenfunctions.

  The case when $\G$ acts freely on $D$ precisely corresponds to the
case discussed in \ref{motivation} (this follows from the
uniformization theorem for the Riemann surface $Y$). Our results
hold for general co-compact subgroup $\G$ (and in fact, with
slight modifications, for any lattice $\G \subset G$).

We will identify the upper half plane $\uH$ (or $D$) with $G / K$,
where $K = PO(2)$ is a maximal compact subgroup of $G$ (this
follows from the fact that $G$ acts transitively on $\uH$ and the
stabilizer in $G$ of the point $z_0 = i \in \uH$ coincides with
$K$).

We denote by $X$ the compact quotient $\G\sm G$ (we call it the
automorphic space). In the case when  $\G$ acts freely  on $\uH$
one can identify the space $X$ with the bundle $S(Y)$ of unit
tangent vectors to the  Riemann surface $Y = \G \sm \uH$.

 The group $G$ acts on $X$ (from the right) and hence on the space of
functions on $X$. We fix the unique $G$-invariant measure $\mu_X$
on $X$ of total mass one. Let $L^2(X)=L^2(X,d\mu_X)$ be the space
of square integrable functions and $(\Pi_X, G, L^2(X))$ the
corresponding unitary representation. We will denote by $P_X$ the
Hermitian form on $L^2(X)$ given by the scalar product. We denote
by $||\ ||_{X}$ or simply $||\ ||$ the corresponding norm and by
$\langle f,g \rangle_X$ the corresponding scalar product.

The identification  $Y=\G\sm \uH\simeq X/K$ induces the embedding
$L^2(Y)\sub L^2(X)$.
 We will always identify the space $L^2(Y)$ with the subspace of
$K$-invariant functions in  $L^2(X)$.

Let  $\phi$ be a normalized eigenfunction of the Laplace-Beltrami
 operator on $Y$. Consider the closed $G$-invariant subspace
$L_\phi\sub L^2(X)$ generated by $\phi$ under the  action of $G$.
It is well-known that $(\pi,L)=(\pi_\phi, L_\phi)$ is an
irreducible unitary representation of $G$ (see \cite{GGPS}).

Usually it is more convenient to work with the space $V = L^\8$ of
smooth vectors in $L$. The unitary Hermitian form $P_X$ on $V$ is
$G$-invariant.

A smooth representation $(\pi, G, V)$ equipped with a positive
$G$-invariant Hermitian form $P$ we will call  a {\it smooth
pre-unitary representation}; this simply means that $V$ is the
space of smooth vectors in the unitary representation obtained
from $V$ by completion with respect to $P$.

   Thus starting with an automorphic function $\phi$ we constructed
   an irreducible smooth pre-unitary representation $(\pi, V)$.
   In fact we constructed this space together with a canonical
   morphism $\nu : V  \to C^\8 (X)$ since $C^\8(X)$ is the smooth
   part of $L^2(X)$.

   \defn {enhanced}A smooth pre-unitary representation
   $(\pi, G, V)$ equipped with a $G$-morphism $\nu: V \to C^\8(X)$
   we will call an {\it $X$-enhanced representation}.

   We will assume that the morphism $\nu$ is normalized,
    i.e. it carries the standard $L^2$ Hermitian
   form $P_X$ on $C^\8(X)$ into Hermitian form $P$ on $V$.

   Thus starting with an automorphic function $\phi$ we
   constructed

   (i) An $X$-enhanced irreducible pre-unitary representation
   $(\pi, V, \nu)$,

   (ii) A $K$-invariant unit vector $e_V \in V$
   (this vector is just our function $\phi$).

   Conversely, suppose we are given an irreducible smooth
   pre-unitary
    $X$-enhanced representation $(\pi,V, \nu)$ of
the group $G$ and a $K$-fixed unit vector $e_V \in V$. Then the
function $\phi = \nu(e_V) \in C^\8(X)$ is $K$-invariant and hence
can be considered as a function on $Y$. The fact that the
representation $(\pi, V)$ is irreducible implies that $\phi$ is an
automorphic function, i.e. an eigenfunction of Laplacian on $Y$.

 Thus we have established a natural correspondence between
 Maass forms $\phi$ and tuples $(\pi, V, \nu,
 e_V)$,   where $(\pi, V,\nu)$ is an $X$-enhanced irreducible
 smooth pre-unitary representation and $e_V \in V$ is a unit
 $K$-invariant vector.

\subsection{Decomposition of the representation $(\Pi_X, G,
L^2(X))$}\label{decompos}

   It is well known that for $X$ compact the
   representation $(\Pi_X, G, L^2(X))$ decomposes into a direct
   (infinite) sum
\begin{equation}\label{spec-L2X}
L^2(X)=\oplus_j (\pi_j, L_j)
\end{equation}
of irreducible unitary representations of $G$ (all representations
appear with finite multiplicities (see \cite{GGPS})). Let $(\pi,
L)$ be one of these irreducible "automorphic" representations and
$V = L^\8$ its smooth part. By definition $V$ is given with a
$G$-equivariant isometric morphism $\nu: V \to C^\8(X)$, i.e. $V$
is an $X$-enhanced representation.

   If $V$ has a $K$-invariant vector it corresponds to a Maass form.
There are  other spaces in this decomposition which  correspond to
discrete series representations. Since they are not related to
Maass forms we will not study them in more detail.

\subsection{Representations of $\PGLR$}\label{irrep}
 All irreducible
unitary representations of $G$ are classified. For simplicity we
consider first those with a nonzero $K$-fixed vector (so called
representations of class one) since only these representations
arise from Maass forms. These are the representations of the
principal and the complementary series and the trivial
representation.

   We will use the following standard explicit model for irreducible
smooth   representations of  $G$.

    For every complex number $\lm$ consider
 the space $V_\lm$ of smooth even homogeneous functions on
$\br^2\sm0$ of homogeneous degree $\lm-1 \ $ (which means that
$f(ax,ay)=|a|^{\lm-1}f(x,y)$ for all $a\in\br \setminus 0$).
 The representation $(\pi_\lm, V_\lm)$ is induced by
the action of the group $\GLR$ given by $\pi_\lm (g)f(x,y)=
f(g\inv (x,y))|\det g|^{(\lm-1)/2}$. This action is trivial on the
center of $\GLR$ and hence defines a representation of $G$. The
representation $(\pi_\lm, V_\lm)$ is called {\it representation of
the generalized principal series}.

When $\lm=it$ is  purely imaginary  the representation $(\pi_\lm
,V_\lm)$ is pre-unitary; the $G$-invariant scalar product in
$V_\lm$ is given by $\langle f,g \rangle_{\pi_\lm}=\frac{1}{2\pi}
\int_{S^1} f\bar g d\te$. These representations are   called
representations of {\it the principal series}.

 When $\lm\in (-1,1)$ the representation $(\pi_\lm ,V_\lm)$ is called
 a representation of the complementary series. These representations
 are also pre-unitary, but the formula for the scalar product is
 more complicated (see \cite{G5}).

 All these representations have $K$-invariant vectors.
 We fix a $K$-invariant unit vector $e_{\lm} \in V_\lm$ to be
  a function which is one on the unit circle in $\br^2$.

  Representations of the principal and the complimentary series exhaust
  all nontrivial irreducible pre-unitary representations of $G$
  of class one. The rest of unitary irreducible representations of
  $G$ could be realized as submodules (or quotients) in the spaces
  $V_\lm$ for $\lm\in\bz$. These are called discrete series
  representations (\cite{G5}, \cite{L}).

 In what follows we will do necessary computations for
 representation of the principal series. Computations for the
 complementary series are a little more involved but essentially
 the same (compare with \cite{BR1}, section 5.5, where similar
 computations are described in detail).

Suppose we are given a class one $X$-enhanced representation
 $\nu: V_{\lm} \to C^\8(X)$; we assume $\nu$ to be an isometric embedding.
 Such $\nu$ gives rise to an
eigenfunction of the Laplacian on the Riemann surface $Y = X/K$ as
before. Namely, if $e_{\lm} \in V_\lm$ is a unit  $K$-fixed vector
then the function $\phi = \nu(e_\lm)$ is a normalized
eigenfunction of the Laplacian on the space $Y = X/K$ with the
eigenvalue $\mu=\frac{1-\lm^2}{4}$. This explains why $\lm$ is a
natural parameter to describe Maass forms. We note that the
Casimir operator $\mathcal{C}$ is a scalar operator on $V_\lm$
with the same eigenvalue. However, eigenspaces of $\mathcal{C}$ in
$C^\8(X)$ correspond only to isotypic components because of
possible multiplicities.

\subsection{Helgason's representation and Frobenius reciprocity}
\label{frob-sect} Here we reformulate Helgason's representation
(\ref{He-int}) for $\G$-invariant eigenfunctions in terms of
Frobenius reciprocity of Gel'fand and Fomin.

Let $(\pi,G, V)$ be an irreducible unitary $X$-enhanced
representation. We have the following {\bf Frobenius reciprocity}
(\cite {GGPS}, \cite{Ol}, \cite{BR2}):

\begin{thm}{Frob-thm}
\begin{align}\label{Frob-iso}
\Mor_G(V, C^\8(X)) \simeq \Mor_{\G}(V,\bc)\ .
\end{align}
\end{thm}
Namely, to every $G$-morphism $\nu: V\to  C^\8(\G\sm G)$
corresponds a $\G$-invariant functional $I$ on the space $V$ given
by $I(v) = \nu(v)(e)$ (here $e$ is the identity in $G$). Given $I$
we can recover $\nu$ as $\nu(v)(g) =I(\pi(g)v)$.

In particular, let $\pi$ be of class one and $e_0\in V$ be a unit
$K$-fixed vector then the corresponding  eigenfunction (i.e. the
Maass form) on $D$ (or $\uH$) is given by
\begin{align}\label{I-Maass}
\phi(z)=I(\pi(g)e_0) \ ,
\end{align}
with $g\cdot 0=z\in D$ (or correspondingly $g\cdot i=z\in\uH$).

This is exactly the Helgason's representation \eqref{He-int} if we
view the automorphic functional $I$ as a distribution on the space
$V\simeq C^\8(S^1)$.

Hence, Helgason's representation shows how to realize the
$K$-fixed vector (i.e. the Maass form) on $D$. However, it does
not show how to realize other vectors in $V$ (and apart from $e_0$
those could not be realized in the space of functions on $D$).
Zelditch \cite{Z2} noticed how to re-write Helgason's
representation in a form appropriate for a general vector $v\in
V$.

Namely, let us choose an identification $V\simeq V_\lm$ and
consider the following left $\G$-invariant distribution on
$D\times B$:
\begin{align}\label{Zel-T}
e_\pi(g)=e^{ (\frac{1+\lm}{2})<z,b>}d{\rm vol}(z)dT(b)
\end{align}
with $g=(z,b)$ under the identification $D\times B\simeq G$ in
\ref{geo-setting}. We have $e_\pi\in\CDD(\G\sm G)$. It is easy to
see that in terms of Frobenius reciprocity this distribution
 is nothing else than
\begin{align}\label{Frob-T}
e_\pi(g)=I(\pi(g)\dl),
\end{align}
where $\dl=\dl_1=\sum_ke_{2k}$ is the distribution which is
formally the sum of all $K$-types in the standard basis of $V$
(see \cite{Z2},\cite{L}) or simply is equal to the Dirac delta
distribution at $1\in S^1$  in the realization $V\simeq
V_\lm\simeq C_{even}^\8(S^1)$. We note (see \cite{L}, \cite{Z2})
that unit vectors $e_{2k}$ become exponents $e_{2k}=\exp(2\pi
i2k\te)$ in the realization $V\simeq C_{even}^\8(S^1)$ of the
principal series representations of $\PGLR$.

Hence, we see that the distribution $e_\pi$ vanishes on functions
which are orthogonal to $V\subset C^\8(X)$ and on $V$ takes value
$1$ on vectors in the standard basis $\{e_{2k}\}$. This is exactly
the description given in \cite{Z2} (Proposition 2.2). We will use
the representation (\ref{Frob-T}) extensively in what follows.

The distribution  $e_\pi$ gives rise to the imbedding
$C_{even}^\8(S^1)\to V\subset C^\8(X)$, $v\mapsto\phi_v(g)$ via
\begin{align}\label{V-imbed}
\phi_v(g)=\int_Ke_\pi(gk)\bar v(k\cdot 1)dk=I(\pi(g)v)\
\end{align}
which again the isomorphism \eqref{Frob-iso}.

\section{Trilinear invariant functionals}\label{3prod-def}

We introduce now the invariant trilinear functionals on
irreducible representations which will be our main tool in what
follows.

\subsection{Automorphic triple products}\label{aut3prod}
Suppose we are given three $X$-enhanced
representations of $G$
\begin{align*}
\nu_j:V_j\to C^\8(X),\ \
 j=1,2,3\ .
\end{align*}

We define the  $G$-invariant trilinear form
$l^{aut}_{\pi_1,\pi_2,\pi_3}:V_1\otimes V_2\otimes V_3\to\bc$ ,
 by formula

\begin{align}\label{aut3prod-def}
l^{aut}_{\pi_1,\pi_2,\pi_3}(v_1\otimes v_2\otimes v_3)= \int_{
X}\phi_{v_1}(x)\phi_{v_2}(x)\phi_{v_3}(x)d\mu_X\ ,
\end{align}
where $\phi_{v_j}=\nu_j(v_j)\in C^\8(X)$ for $v_j\in V_j$.

\subsection{Uniqueness of triple products} The central fact about
invariant trilinear functionals is the following uniqueness
result:

\begin{thm}{ubi} Let $(\pi_j,V_j),\ \ j=1,2,3\ ,$
be three irreducible smooth admissible
  representations of  $G$. Then
$\dim\Hom_G(V_1\otimes V_2\otimes V_3,\bc)\leq 1$.
\end{thm}

\begin{rem}{rem-uni}
The  uniqueness statement was proven by Oksak in \cite{O} for the
group $\SLC$ and the proof could be adopted for $\PGLR$ as well
(see also \cite{Mo} and \cite{Lo} for different proofs). For the
$p$-adic $GL(2)$ more refined results were obtained by Prasad (see
\cite{Pr}). He also proved the uniqueness when at least one
representation is a discrete series representation of $\GLR$.

 There is no uniqueness of trilinear functionals for representations
of $\SLR$ (the space is two-dimensional). This is the reason why
we prefer to work with $\PGLR$ (although the method could be
easily adopted to $\SLR$).

For $\SLR$ one has the following uniqueness statement instead. Let
$(\pi,V)$ and $(\sigma, W)$ be two  irreducible smooth pre-unitary
representations of  $\SLR$ of class one. Then the space of
$\SLR$-invariant trilinear functionals on $V\otimes V\otimes W$
which are symmetric in the first two variables is one-dimensional.
This is the correct uniqueness result needed if one wants to work
with $\SLR$.
\end{rem}

\subsection{Model trilinear functionals} \label{modfunc}

For every $\lm \in \bc$ we denote by $(\pi_\lm,V_\lm)$  the smooth
class one representation of the generalized principle series of
the group $G=\PGLR$ described in \ref{irrep}. We will use the
realization of $(\pi_\lm, V_\lm)$ in the space of smooth
homogeneous functions on $\br^2 \setminus 0$ of homogeneous degree
$\lm - 1$ .

For explicit computations it is often convenient to pass from
plane model to a circle model. Namely, the  restriction of
functions in $V_\lm$ to the unit circle $S^1 \subset \br^2$
defines an isomorphism of the space $V_\lm$ with the space
$C^\8(S^1)^{even}$ of even smooth functions on $S^1$ so we can
think about vectors in $V_\lm$ as functions on $S^1$.

We describe now the {\it model} invariant trilinear functional
using the explicit geometric models for irreducible
representations described above. Namely, for given three complex
numbers $\lm_j$, $j=1, 2, 3$, we construct explicitly nontrivial
trilinear functional $\  l^{mod}: V_{\lm_1} \otimes V_{\lm_2}
\otimes V_{\lm_3} \to \bc$ by means of its kernel.

\subsubsection{Kernel of  $l^{mod}$ }\label{kernel}

 Let
$\om(\xi,\eta)= \xi_1\eta_2-\xi_2\eta_1$ be $\SLR$-invariant of a
pair of
 vectors $\xi,\ \eta\in\br^2$. We set

\begin{equation}\label{kern}
K_{\lm_1,\lm_2,\lm_3}(s_1,s_2,s_3)= |\om(s_2,s_3)|^{(\al-1)/2}\
 |\om(s_1,s_3)|^{(\beta-1)/2}|\om(s_1,s_2)|^{(\g-1)/2}
\end{equation}
for $s_1,s_2,s_3 \in\br^2\sm 0$, where $ \al =\lm_1-\lm_2-\lm_3,
  \beta= -\lm_1+\lm_2-\lm_3,\g = -\lm_1-\lm_2+\lm_3 $.

 The kernel function $K_{\lm_1,\lm_2,\lm_3}(s_1,s_2,s_3)$
satisfies two  main properties:
\begin{enumerate}
\item $K$ is invariant with respect to the diagonal action of
$\SLR$.

\item $K$ is homogeneous of degree $-1 -\lm_j$ in each variable
$s_j$.
\end{enumerate}

 Hence
if $f_j$ are homogeneous functions of degree $-1+\lm_j$, then the
function
$$F(s_1,s_2,s_3)=
f_1(s_1)f_2(s_2)f_3(s_3)K_{\lm_1,\lm_2,\lm_3}(s_1,s_2,s_3)\ ,$$ is
homogeneous of degree $-2$ in each variable $s_j\in\br^2\sm 0$.

\subsubsection{Functional $l^{mod}$ }\label{int-of-kern}
To define the model trilinear functional $l^{mod}$ we notice that
on the space $\CV$ of functions of homogeneous  degree $-2$ on
$\br^2\sm 0$ there
 exists a natural $\SLR$-invariant
functional $\fL:\CV \to\bc$ . It is given by the formula $\fL(f) =
\int_\Sigma f d\sigma$ where the integral is taken over any closed
curve $\Sigma\subset \br^2\sm 0$ which goes around $0$  and the
measure $d\s$  on $\Sigma$ is given by the area element inside of
$\Sigma$ divided by $ \pi$; this last normalization factor is
chosen    so that $\fL(Q^{-1}) = 1$ for the standard quadratic
form $Q$ on $\br^2$.

Applying $\fL$ separately to each variable $s_i\in \br^2\sm 0$ of
the function $F(s_1,s_2,s_3)$ above we obtain the $G$-invariant
functional
\begin{equation}\label{pairing}
l^{mod}_{\pi_1,\pi_2,\pi_3}(f_1\otimes f_2\otimes f_3):= \langle
\fL\otimes\fL\otimes\fL ,F \rangle\ .
\end{equation}
We call it the {\it model triple product} and denote  by
$l^{mod}_{\pi_1,\pi_2,\pi_3}$.

In the circle model this functional is expressed by the following
integral:
\begin{equation}\label{circleintegral}
l^{mod}_{\pi_1,\pi_2,\pi_3}(f_1\otimes f_2\otimes f_3)= (2
\pi)^{-3} \iiint f_1(x)f_2(y)f_3(z)K_{\lm_1,\lm_2,\lm_3} (x,y,z)
dx dy dz ,
\end{equation}

where $x, y , z\in S^1$ are the standard angular parameters on the
circle and
\begin{equation}\label{circle-kernel}
 K_{\lm_1,\lm_2,\lm_3} (x,y,z) =|\sin(y-z)|^{(\al-1)/2}
     |\sin(x-z)|^{(\beta-1)/2} |\sin(x-y)|^{(\g-1)/2}\
\end{equation}
with $\al,\ \beta,\ \g \in i\br$ as before.

\rem {remark} The integral defining the trilinear functional is
often divergent and the functional should be defined using
regularization of this integral. There are standard procedures how
to make such a regularization (see e.g. \cite{G1}).

\subsection{Coefficients of proportionality}\label{coef-prop} By the uniqueness
principle, for automorphic representations $\pi_1,\pi_2,\pi_3$
there exists a constant $a_{\pi_1,\pi_2,\pi_3}$ of proportionality
between the {\it model} functional (\ref{pairing}) and the {\it
automorphic} functional (\ref{aut3prod-def}) :
\begin{align}\label{coef-a-def}
l^{aut}_{\pi_1,\pi_2,\pi_3}=a_{\pi_1,\pi_2,\pi_3} \cdot
l^{mod}_{\pi_1,\pi_2,\pi_3}\ .
\end{align}

\subsubsection{Bounds on $a_{\pi_1,\pi_2,\pi_3}$}\label{a-bound}
 In this paper we will need the following particular
case of a general problem of estimating the coefficients
$a_{\pi_1,\pi_2,\pi_3}$. Let us fix an automorphic representation
$\pi_1\simeq \pi_\mu$ and let $\pi_2=\pi_3\simeq \pi_{\lm_i}$ as
$|\lm_i|\to\8$ through the set of parameters of automorphic
representations of class one. We have the following (the so-called
convexity) bound:

\begin{prop}{propklambda}  There exists an effective
constant $C$ such that for any $\pi_\mu$ and $\pi_{\lm_i}$
\begin{align}\label{coef-a-bound}
|a_{\pi_\mu,\pi_{\lm_i},\pi_{\lm_i}}| \leq
C(\max(|\mu|,|\lm_i|))^\haf \ .
\end{align}
\end{prop}

\begin{proof} This follows from methods of \cite{BR3}.
For $\mu$ fixed and $|\lm_i|\to\8$ this is also shown in
\cite{Re2} by a slightly different argument. We discuss similar
bounds for the case of the representation $\pi_\mu$ of discrete
series in the course of the proof of Theorem \ref{invar-thm}.
~\end{proof}

\subsubsection{A conjecture} \label{a-conj}

The major problem in the theory of automorphic functions and
analysis on $Y$ is to find a method which would allow one to
obtain better bounds for coefficients $a_{\pi_1,\pi_2,\pi_3}$.

We would like to make the following conjecture concerning the size
of coefficients $a_{\pi_\mu,\pi_{\lm_i},\pi_{\lm_i}}$:

\begin{conj}{CONJ}
 For a fixed $\pi_\mu$ and for any
$\eps>0$ there
exists $C_\eps>0$ independent of $\lm_i$ such that
$$|a_{\pi_\mu,\pi_{\lm_i},\pi_{\lm_i}}|\leq C_\eps|\lm_i|^\eps\ ,$$
as $|\lm_i|\to\8$.
\end{conj}

In a special case of a congruence subgroup $\G$ this conjecture is
consistent with the Lindel\"{o}f conjecture from the theory of
automorphic $L$-functions (see \cite{Wa} for a connection to the
theory of $L$-functions, \cite{Sa2} for the survey and \cite{BR1},
\cite{BR3}, \cite{Re2} for the connection to trilinear
functionals).

\section{PDO and microlocal distributions $dU_i$}\label{pdo-3lin-sect}

In this section we piece together pseudo-differential operators
and representation theory in order to express Zelditch's
microlocal lifts of eigenfunctions in terms of representation
theory.

\subsection{PDO} Let $\phi\in C^\8(Y)$ be an eigenfunction with
the eigenvalue $\mu=\frac{1-\lm^2}{4}$, $(\pi,V)$ the
corresponding automorphic representation with the automorphic
functional $I\in V^*$ and the Helgason-Zelditch distribution
$e_\pi\in \CDD(X)$ (see (\ref{Zel-T})). Let also $a(g)=a(z,b)\in
C^\8(\G\sm D\times B)$ be a symbol of order zero (which we assume
is independent of $\lm$).

In \cite{Z2}, on the basis of the representation
(\ref{Op-inv-def}), Zelditch defined the corresponding
pseudo-differential operator $A=Op(a):C^\8(Y)\to C^\8(Y)$ by
\begin{align}\label{pdo-oper-Y}
Op(a)\phi(z)=\int_B a(z,b)e^{ (\frac{1+\lm}{2})<z,b>}dT(b)\ .
\end{align}
We can rewrite this in the form
\begin{align}\label{pdo=mult}
Op(a)\phi(z)=\int_K a(gk)e_\pi(gk)dk\ .
\end{align}
Hence, the action of $A$ on $\phi$ reduces to the multiplication
of the corresponding distribution $e_\pi$ by the symbol $a(g)$ and
then taking the $K$-invariant part of the result.

\subsection{Distributions $dU_i$} Interpreting pseudo-differential
operator as an observable in Quantum Mechanics one is led to the
introduction of correlation functions or matrix coefficients.
Namely, one is interested in studying following quantities
\begin{align}\label{matrix-coef}
A_{ij}=<Op(a)\phi_i,\phi_j>\ .
\end{align}
One view these as distributions on the space of symbols. We will
concentrate on the diagonal terms $<Op(a)\phi_i,\phi_i>$ first.
This leads us to the following definition of distributions $dU_i$
on $X$ associated to eigenfunctions $\phi_i$ on $Y$:
\begin{align}\label{dU}
<Op(a)\phi_i,\phi_i>:=\int_X a(x)dU_i\ .
\end{align}
Using the interpretation (\ref{pdo=mult}) we arrive to the
following defining relation for the distributions $dU_i$:
\begin{align}\label{dU-X}
\int_Xa(x)e_{\pi_i}(x)\bar\phi_i(x)dx:=\int_X a(x)dU_i\ .
\end{align}
Hence we see that
\begin{align}\label{dU=e*f}
dU_i=e_{\pi_i}(x)\bar\phi_i(x)
\end{align}
as distributions on $X$. Note that from the construction of
$e_{\pi_i}$ it follows that $\int_X1dU_i=\int_X|\phi_i|^2dx=1$.

\subsection{Automorphic functions as symbols}\label{main-proof}
We now rephrase  Theorem \ref{main-thm} from the Introduction in
terms of automorphic representations (while in the Introduction we
stated it in equivalent terms of eigenspaces of the Casimir; see
Section \ref{reps}). This theorem underlies our  study of action
of pseudo-differential operators on eigenfunctions.

\begin{thm}{main-thm2} Let $\nu_\mu:C^\8(S^1)\to V_\mu \subset C^\8(X)$ be an
irreducible automorphic representation and the corresponding
$G$-morphism. For any $a\in V_\mu$ and any Maass forms $\phi_i$
and $\phi_j$ there are a constant $a_{\mu,\mu_i,\mu_j}$ and an
explicit distribution $l_{\mu,\mu_i,\mu_j}\in \mathcal{D}(S^1)$,
depending only on ${\mu,\mu_i,\mu_j}$ but not on the choice of the
corresponding functions, such that the following relation holds
\begin{align}\label{Op-tripl2}
<Op(a)\phi_i,\phi_j>_Y=a_{\mu,\mu_i,\mu_j}\cdot
<l_{\mu,\mu_i,\mu_j},v_a>_{V_\mu} ,
\end{align}
where $a=\nu_\mu(v_a)$ with $v_a\in C^\8(S^1)$.
\end{thm}
\begin{proof} We deal with symbols coming from the class one representations
$V_\mu$. The case of symbols coming from discrete series is
similar and is explained in detail in the course of the proof of
Theorem \ref{invar-thm}.

  As the symbol $a\in V_\mu$ belongs to an irreducible
representation we have from \eqref{pdo=mult} and
\eqref{aut3prod-def}
\begin{align}\label{Op-triple2}
<Op(a)\phi_i,\phi_j>=\int_X
\psi(x)e_{\pi_i}(x)\bar\phi_j(x)dx=l_{\pi_\mu,\pi_i,\pi_j}^{aut}(a,e_{\pi_i},
\phi_j)\ .
\end{align}
Hence from \eqref{Frob-T} and \eqref{coef-a-def} we have
\begin{align}\label{dist-main-thm}
l_{\pi_\mu,\pi_i,\pi_j}^{aut}(a,e_{\pi_i}, \phi_j)=\
a_{\pi_\mu,\pi_{\lm_i},\pi_{\lm_i}}\cdot
l_{\pi_\mu,\pi_i,\pi_j}^{mod}(v_a,\dl, e_0)\
\end{align}
and hence $l_{\pi_\mu,\pi_i,\pi_j}^{mod}(v_a,\dl,
e_0)=<l_{\mu,\mu_i,\mu_j},v_a>_{V_\mu} $ could be viewed as the
evaluation of the distribution with the explicit kernel on $S^1$,
given by \eqref{integ-reduc2}, on the function $v_a$ once we
consider the identification $C^\8(S^1)_{even}\simeq V_\mu$.

\end{proof}

\section{Invariance of $dU_i$ under the geodesic flow}

\subsection{Geodesic flow} It is well known
that under the identification $S^*(Y)\simeq X$ the geodesic flow
$G_t$ on $S^*(Y)$ corresponds to the right action on $X$ of the
diagonal subgroup $T=\{g_t=diag(e^t,e^{-t})|\ t\in \br\}\subset G$
(see \cite{GF}).

\subsection{Asymptotic invariance}
 In order to prove
asymptotic invariance of distributions $dU_i$ we will show that
for any symbol $a(x)\in C^\8(X)$
\begin{align}\label{invarint-general}
\left|\int_X\left(a(x)-a(xg_t)\right)dU_i\right|=O_{a,t}(|\lm_i|^{-\al})\
,
\end{align}
for some $\al>0$ and with a uniform constant in the $O$-term as
$t$ changing in a compact set and $a$ bounded (w.r.t. a Sobolev
norm on $C^\8(X)$). We will show that one can choose $\al=1$
above. Such bounds are usually obtained as a consequence of the
Egorov-type theorem (see \cite{Z1}). Zelditch found another way to
prove such bounds based on the exact differential relation
satisfied by $dU_i$. We will use trilinear invariant functionals
introduced above in order to prove (\ref{invarint-general}). As a
consequence of our proof we will be able to speculate (on the
basis of Conjecture \ref{CONJ}) that the true order of decay in
(\ref{invarint-general}) should be $|\lm_i|^{-3/2}$. We note that
the Egorov's theorem gives only $|\lm_i|\inv$ as the order of
decay in (\ref{invarint-general}).

In order to be able to connect distributions $dU_i$ to trilinear
invariant functionals we will consider symbols which are
themselves automorphic functions (i.e. symbols which belong to one
of automorphic representations $V_i$). Such functions  are dense
in $C(X)$ (e.g. union of basis $\{\psi_k^i(x)\}$ of all spaces
$V_i$) and hence describe $dU_i$ uniquely.

We have the following

\begin{thm}{invar-thm} For any fixed automorphic
representation $(\pi_\mu, V_\mu)$ there exists an explicit
constant $c_\mu$ such that for a given automorphic function (which
we view as a symbol) $\psi(x)\in V_\mu\subset C^\8(X)$ the
following relation holds
\begin{align}\label{invar-bound}
<Op(\psi)\phi_i,\phi_i>=a_{\pi_\mu,\pi_{\lm_i},\pi_{\lm_i}}|\lm_i|^{-\haf}
c_\mu d_\mu(\psi)\ +\ O_{\psi,\mu}(|\lm_i|\inv)\ ,
\end{align}
where $d_\mu$ is the properly normalized, independent of
 $\lm_i$, $T$-invariant functional on $V_\mu$.

The constant in the $O$-term above is effective in $\mu$ and (the
Sobolev norm of) $\psi$.
\end{thm}

\begin{cor}{invar-cor}For $\psi$ as above we have
\begin{align}\label{invariant-auto}
\left|\int_X\left(\psi(x)-\psi(xg_t)\right)dU_i\right|=O_{\psi,t}(|\lm_i|\inv)\
,
\end{align}
as $|\lm_i|\to\8$.
\end{cor}

\begin{rem}{Remark} 1. From the proof it follows that the constant in
the $O$-terms above is bounded by the second $L^2$-Sobolev norm of
$\psi$ and $|\mu|^\haf$. Hence we have for a general symbol $a\in
C^\8(X)$:
\begin{align}\label{invariant-a}
\left|\int_X\left(a(x)-a(xg_t)\right)dU_i\right|\leq C_t\cdot
S_2(a)\cdot |\lm_i|\inv\ ,
\end{align}
as $|\lm_i|\to\8$. Here $S_2$ is the second $L_2$-Sobolev norm on
$X$. This should be viewed as a quantitative version of the
asymptotic invariance of the distributions $dU_i$.

2. We have seen in \eqref{coef-a-bound} that
$|a_{\pi_\mu,\pi_{\lm_i},\pi_{\lm_i}}|\leq C|\lm_i|^{\haf}$ and
hence the coefficients in front of the distribution $d_\mu$ are
uniformly bounded in $\lm_i$. One expects that coefficients
$a_{\pi_\mu,\pi_{\lm_i},\pi_{\lm_i}}$ grow at a much slower rate
(e.g. Conjecture \ref{CONJ}). This is known as the (effective)
Quantum Unique Ergodicity conjecture of Rudnick and Sarnak
solution to which (in an ineffective form) was recently announced
by E. Lindenstrauss \cite{Li2}.
\end{rem}

\subsection{Proof} In order to prove (\ref{invar-bound}) we use our
interpretation of pseudo-differential operators with automorphic
symbols in terms of trilinear functionals.  We first deal with
symbols $\psi$ coming from automorphic representations of class
one described in \ref{irrep}. We have from (\ref{dU}) and
(\ref{dU=e*f})
\begin{align}\label{Op-triple}
<Op(\psi)\phi_i,\phi_i>=\int_X
\psi(x)e_{\pi_i}(x)\bar\phi_xdx=l_{\pi_\mu,\pi_i,\pi_i}^{aut}(\psi,e_{\pi_i},
\phi_i)\ .
\end{align}

Let $\psi=\nu_\mu(v)$ for $v\in V_\mu\simeq C^\8_{even}$, $\dl\in
V_{\lm_i}^*$ the distribution which corresponds to $e_{\pi_i}$ and
$e_0=e_{0,\lm_i}$ the $K$-fixed vector in $V_{\lm_i}$. From the
uniqueness of trilinear functionals (\ref{coef-a-def}) we arrive
at
\begin{align}\label{Op-triple-mod}
<Op(a)\phi_i,\phi_i>=a_{\pi_\mu,\pi_{\lm_i},\pi_{\lm_i}}\cdot
l_{\pi_\mu,\pi_i,\pi_i}^{mod}(v,\dl, e_0)\ .
\end{align}

We use now the explicit description of $l^{mod}$ in
(\ref{circleintegral}) in order to compute the right hand part of
(\ref{Op-triple-mod}). We have (recall that $e_0$ is the constant
function)
\begin{align}\label{integ-reduc1}
 l^{mod}_{\pi_\mu,\pi_{\lm_i},\pi_{\lm_i}}(v,\dl, e_0)=&&\\
 (2\pi)^{-3} \int_{(S^1)^3} v(x)\dl(y)&&
 |\sin(y-z)|^{(\al-1)/2}|\sin(x-z)|^{(\beta-1)/2}
 |\sin(x-y)|^{(\g-1)/2}\ dx dy dz ,\nonumber
\end{align}

where $x, y , z\in S^1$ are the standard angular parameters on the
circle and $\al=-2\lm_i+\mu,\ \beta=-\mu,\ \g=-\mu $ as before .

We note that $\dl=\dl_0$ is the Dirac delta at $y=0$ and hence we
need to compute the following integral
\begin{align}\label{integ-reduc2}
 \int_{(S^1)^2} v(x)|\sin(z)|^{-\lm_i+\haf\mu-\haf}|\sin(x-z)|^{-\haf\mu-\haf}
 |\sin(x)|^{-\haf\mu-\haf}\ dx dz\ .
\end{align}
as $|\lm|=|\lm_i|\to\8$. We compute it with the help of the
stationary phase method. We write (\ref{integ-reduc2}) as
\begin{align}\label{integ-reduc3}
 \int_{(S^1)^2} A(x,z)e^{i\lm p(z)} dx dz\ ,
\end{align}
with the amplitude
\begin{align}\label{amplit}
 A(x,z)=v(x)|\sin(z)|^{\haf\mu-\haf}|\sin(x-z)|^{-\haf\mu-\haf}
 |\sin(x)|^{-\haf\mu-\haf}\
\end{align}
and the phase (which depends only on $z$)
\begin{align}\label{phase}
 p(z)=-\ln|\sin(z)|\ .
\end{align}
 A direct computation shows that the phase $p$ has two
 non-degenerate critical points $z_{\pm}=\pm\pi/2$ with
 the equal contribution to the integral because of the symmetry (all
  functions involved are even on $S^1$). We note now that this
explains why the integral (or more precisely its leading term)
(\ref{integ-reduc1}) gives rise to a distribution on $S^1\times
S^1$ which is invariant under the diagonal action of the diagonal
subgroup $T$. The simple geometric reason for this is that by the
stationary phase method the leading term is given by the value at
the critical points $y=0,\ z=\pi/2$ and $y=0,\ z=-\pi/2$. These
points are fixed points of the diagonal action of $T$. Near these
fixed points elements of $T$ contract (or expand) in direction of
$y$ and expand (respectively contract) by the same amount in
direction of $z$. Hence, distribution supported in these fixed
points is invariant with respect to the diagonal action in the
space $V_{\lm_i}\times V_{\lm_i}\simeq C^\8(S^1\times S^1)$.

On the formal level, the second derivative of the phase at the
fixed point is $1$ and the value of the phase is $0$ and hence the
(equal) contribution from critical points to the integral
(\ref{integ-reduc3}) is given by $|\lm|^{-\haf}$ times the value
of
\begin{align}\label{contrib}
 \int_{S^1} v(x)|\sin(x)|^{-\haf\mu-\haf}|\sin(x-\pi/2)|
 ^{-\haf\mu-\haf}dx=\int_{S^1} v(x)|\sin(2x)|^{-\haf\mu-\haf} dx\ ,
\end{align}
which is exactly a $T$-invariant distribution  on $V_\mu$. Let
$d_\mu$ be the unique $T$-invariant distribution on $V_\mu$
normalized to have value one on the $K$-fixed vector $e_0$. The
value of the  functional in \eqref{contrib} on $e_0$ (i.e. on the
constant function $1$) is given by the classical integral and
gives the main term in \eqref{invar-bound}
\begin{align}\label{k-val}
 c_\mu=\int_{S^1} |\sin(2x)|^{-\haf\mu-\haf} dx\ = \
 \frac{2^{-\haf+\haf\mu}\G(\haf-\haf\mu)}{\G(\frac{3}{4}-\qtr\mu)^2}\ .
\end{align}
We note that all but finite number of eigenvalues of $\Dl$ are
greater than $\qtr$ and hence all but finite number of $\mu$'s are
purely imaginary. From the Stirling's formula we have
$|c_\mu|=O(|\mu|^\haf)$ as $|\mu|\to\8$.

The reminder in the stationary phase method is of order
$O_v(|\lm|^{-3/2})$. Here the constant in the $O$-term is bounded
by the first derivative of $v$. From  \eqref{coef-a-bound} we see
that $|a_{\pi_\mu,\pi_{\lm_i},\pi_{\lm_i}}|\ll |\lm_i|^{\haf}$ and
hence the $O$-term claimed.

We now turn to symbols coming from the discrete series
representations. Let $k>1$ be an odd positive integer and $D_k$
the space (of smooth vectors) of the corresponding discrete series
representation (see \cite{G5} and \cite{L} for various
descriptions of discrete series). We will use the following
well-known realization of discrete series. Let $V_{-k}$ be the
space of smooth functions of the homogeneous degree $-k-1$ on
$\br^2$. The space $D_k$ could be realized as a subspace in
$V_{-k}$ (with the quotient isomorphic to the finite-dimensional
representation of the dimension $k$; see \cite{G5}). We denote the
corresponding imbedding by $i_k:D_k\to V_{-k}$.

For a representation $V_{\lm_i}$ of the principal series we need
to construct a (unique) $G$-invariant functional on the tensor
product $D_k\otimes V_{\lm_i}\otimes V_{\lm_i}$ . We first note
that the formula \eqref{kern} defines the kernel of the
$G$-invariant functional on the (reducible) representation
$V_{-k}\otimes V_{\lm_i}\otimes V_{\lm_i}$. For this one can use
general methods of analytic continuation of integrals described in
\cite{G1} to regularize the integral in \eqref{circleintegral}.
This gives meaning to the value of this integral for $\mu=-k$.
This is true for any $\lm$ which is not a pole of the analytic
continuation of \eqref{circleintegral}. For a value of $\lm$ which
is a pole of \eqref{circleintegral} one also can assign an
invariant functional by taking the residue. It is easy to see that
$\mu=-k$ is not a pole for the analytic continuation. This is
especially easy to see for $\lm$ non-real. Hence we have an
invariant functional on $V_{-k}\otimes V_{\lm_i}\otimes V_{\lm_i}$
which gives rise to the invariant functional on the subspace
$D_k\otimes V_{\lm_i}\otimes V_{\lm_i}$. It is easy to see that
such a functional is non-zero on $D_k$ for $k\equiv 3$ (mod $4$).
For $k\equiv 1$ (mod $4$) this functional vanishes on $D_k$ and
one have to consider the derivative of $l_{\mu,\lm_i,\lm_i}$ in
$\mu$ evaluated at $\mu=-k$. In both cases the functionals
obtained are very similar. We consider the case $k\equiv 3$ (mod
$4$) for simplicity and leave the similar case $k\equiv 1$ (mod
$4$) to the reader. The value of the described above invariant
functional on the triple $v\in D_k$, $\dl,\ e_0\in V_{\lm_i}$ is
defined by the analytic continuation of the integral
\begin{align}\label{integ-D-ser}
 \int_{(S^1)^2} v(x)|\sin(z)|^{-\lm_i-\haf k-\haf}|\sin(x-z)|^{\haf k-\haf}
 |\sin(x)|^{\haf k-\haf}\ dx dz\ .
\end{align}
The value of \eqref{integ-D-ser} is obtained by the analytic
continuation of the distribution $f_s$ on $S^1$ which is given by
the kernel $|\sin(z)|^s$ for $Re(s)>-1$ and then analytically
continued to $s={-\lm_i-\haf k-\haf}$. Moreover, the contribution
from a small neighborhood of singular points ($z=0,\ \pi$) to the
value of $f_s$ on any fixed smooth function is negligible as
$Im(s)\to\8$. Namely, let $g(z)\in C^\8(S^1)$ be a smooth function
with a support in $(-0.1\pi,0.1\pi)$ then $|f_s(g)|\ll
|Im(s)|^{-N}$ for any $N>0$ as $Im(s)\to\8$. This implies that we
can disregard any small enough, fixed neighborhood of $z=0$ in the
integral \eqref{integ-D-ser} and hence we end up with the integral
without non-integrable singularities. Such an integral could be
treated in the same fashion as before and hence the leading term
is given by $|\lm_i|^{-\haf}$ times the value of the integral
\begin{align}\label{contrib-d-ser}
 \int_{S^1} v(x)|\sin(2x)|^{\haf k-\haf} dx\ .
\end{align}
This is again the unique (up to a constant) $T$-invariant
distribution on $D_k$. Taking into account that the reminder in
the stationary phase method is of order $|\lm|^{-3/2}$ we arrive
to the reminder in \eqref{invariant-auto} for a fixed $k$. We
study now the dependence on $k$ of the constant in the reminder.
For this we normalize this $T$-invariant distribution by computing
its value on a special vector in $D_k$. Namely, let
$w_k=\exp(i(-1-k))$ be the highest weight vector in $D_k$
(strictly speaking w.r.t. $\PSLR$). Let $d_k$ be the distribution
taking the value $1$ on a unit vector proportional to $w_k$. The
value of \eqref{contrib-d-ser} on $w_k$ is given by the classical
integral (due to Ramanujan, see \cite{Ma})
\begin{align}\label{v-k-val}
\int_{S^1} |\sin(2x)|^{\haf k-\haf}e^{i(-1-k)x} dx\ = \ e^{-\qtr
i\pi (k+1)} \frac{2^{\haf-\haf k}\G(\haf+\haf k)}{\G(1+\haf
k)\G(\haf)} \ .
\end{align}
From this we see that the last expression is of order
$\al_k=\pi^{-\haf}2^{\haf-\haf k}|k|^{-\haf}$. Taking into account
that $||w_k||^2_{D_k}=\G(2k)\inv$ (see \cite{G5}) we arrive at
$|c_k|=\pi^{-\haf}2^{\haf-\haf k}|k|^{-\haf}\G(2k)^\haf$.

To estimate the reminder we need to estimate the automorphic
coefficient $a_{\pi_k,\pi_{\lm_i},\pi_{\lm_i}}$. We will show that
the bound
\begin{align}\label{discr-ser-a}
|a_{\pi_k,\pi_{\lm_i},\pi_{\lm_i}}|^2=O(|k|2^k\G(2k)\inv)
\end{align}
 holds. This bound is similar to the bound
\eqref{coef-a-bound}. The appearance of the $\G$-function is due
to the awkward normalization of the trilinear functional for the
discrete series. This is mostly due to the author's lack of
knowledge of good models of discrete series. We expect that a
stronger bound follows from methods of \cite{BR3}. We show here
how to obtain the claimed bound by a more elementary means.

Let $k\ll |\lm|$. We estimate the value of the model trilinear
invariant functional $l^{mod}_{\pi_k,\pi_{\lm},\pi_{\lm}}$ on
specially chosen  (smooth) vectors. For the automorphic trilinear
functional we use the bound coming from the maximum modulus
estimate on vectors in the automorphic representation $D_k$. This
will give us a bound on the coefficient of proportionality
$a_{\pi_k,\pi_{\lm_i},\pi_{\lm_i}}$.

We choose the triple $w_k\otimes e_0\otimes e_{k+1}$, where $w_k$
is as above and $e_l\in V_{\lm_i}$ is a unit vector of the
$K$-type $l$ which we will view as a function $e_l=\exp(il\te)$ in
the realization $V_{\lm_i}\simeq C^\8(S^1)$. We note that since
this triple is invariant under the action of the diagonal copy of
$K$ the integral we have to compute could be reduced to
\begin{align}\label{integ-D-ser2}
 \int_{(S^1)^2} e_0(y)e_{k+1}(z)|\sin(y-z)|^{-\lm_i-\haf k-\haf}
 |\sin(y)|^{\haf k-\haf}
 |\sin(z)|^{\haf k-\haf}\ dy dz\ .
\end{align}
As before the stationary phase method imply that the main
contribution to this integral is given by $|\lm|^{-\haf}$ times
the value of the integral along the line $x-y=\pi/2$, namely
\begin{align}\label{contrib-d-ser2}
 \int_{S^1} |\sin(2x)|^{\haf k-\haf}e^{i(k+1)x} dx\
\end{align}
which we computed above and saw that it is of order of
$\al_k=\pi^{-\haf}2^{\haf-\haf k}|k|^{-\haf}$. As we mentioned the
norm of $w_k$ is equal to $\G(2k)^{-\haf}$ and hence using Sobolev
type bound from \cite{BR2} we arrive at the pointwise bound for
the automorphic realization $\phi_k(g)=\nu_k(w_k)(g)$ of the
highest weight vector $w_k$ in the discrete series $D_k$ of the
type $\sup_X|\phi_k|\leq C|k|^\haf\G(2k)^{-\haf}$ and hence the
bound on $a_{\pi_k,\pi_{\lm_i},\pi_{\lm_i}}$ claimed in
\eqref{discr-ser-a}. Combined with the computed value for $c_k$,
this gives the bound for the constant in the reminder in
\eqref{invariant-auto}.

\rightline{$\Box$}

\section{Non-negative microlocal lifts}
We now want to correct distributions $dU_i$ by a smaller order
term in $\lm_i$ (as $|\lm_i|\to\8$) in order to obtain probability
measures $dm_i$ on $X$. Namely, we want to construct a family of
probability measures $dm_i$ such that for any $f\in C^\8(X)$ the
following relation holds
\begin{align}\label{dm-dU}
 \int_Xf(x)dU_i=\int_Xf(x)dm_i\ +\ O_f(|\lm_i|^{-\haf})\ .
\end{align}
This is usually done by means of averaging over a small set in the
phase space. Such a procedure is  called Friedrichs symmetrization
(\cite{Sh},\cite{CdV},\cite{Z2}). However, Friedrichs
symmetrization does not commute with the action of $G$ and hence
does not preserve automorphic representations (this problem is
discussed in \cite{Z3}). In this section we show that one can
exhibit a variety of families of probability measures which are
asymptotic to $dU_i$ and constructed via representation theory.

\subsection{Probability measures}
We construct asymptotic to $dU_i$ probability measures on $X$ by
taking restrictions of automorphic functions
$\psi\otimes\bar\psi\in V_{\lm_i}\otimes V_{\lm_i}$ on $X\times X$
to the diagonal $\Dl X\hookrightarrow X\times X$. Where $\psi\in
V_{\lm_i}$ is a specially chosen $L^2$-normalized automorphic
functions. This will give rise to a probability measures since
representations $V_{\lm_i}$ are self-dual and hence the resulting
function is non-negative on $\Dl X$. Our construction is motivated
by Wolpert's approach to the microlocal lift via the Fej\'{e}r
kernel (\cite{Wo}).

Let $\chi(t)\in C^\8(S^1)$ be a smooth non-negative function
supported in $(-\pi/4,\pi/4)$ and with the norm
$\int|\chi(t)|^2dt=1$. We consider a family of vectors $v_r\in
V_{\lm_i}\simeq C^\8_{even}(S^1)\simeq C^\8(S^1)$ for $r>1$
defined by $v_r=2^{-\haf}r^\haf (\chi(rt)+\chi(r(t-\pi/2))$, $t\in
S^1$ (i.e. the sum of two contracted bump functions around $0$ and
around $\pi/2$). We note that $||v_r||=1$. Clearly the function
$\rho_r(x)=\rho_r^{\lm_i}(x)=\nu(v_r)\otimes\nu(v_r)\mid_{\Dl X}$
is a density of a probability measure on $X$. We have the
following

\begin{thm}{asym-dm} For any (symbol) $\psi\in C^\8(X)$ and any
$\eps>0$  there exists an effective constant $C=C_{\psi,\eps}$
such that for any $\lm_i$ we have
\begin{align}\label{aprox-dU-dm}
\left|\int_X\psi(x)dU_i\ - \ \int_X\psi(x)\rho_r(x)dx\right|\ \leq
C|\lm_i|^{-\haf}\ ,
\end{align}
as $r\to\8$ and $|\lm_i|\to\8$ condition to $r\leq
|\lm_i|^{\haf-\eps}$.
\end{thm}
\begin{proof}
Consider a given $r>1$. We may assume that $\psi\in V_\mu$ and
$\psi=\nu(v)$ as in \ref{enhanced}. The integration in
\eqref{invar-bound} is over the small neighborhood (depending on
the value of $r$) of four points $x=0$ or $\haf\pi$ and $y=0$ or
$\haf\pi$. However, as we saw in the proof of Theorem
\ref{invar-thm}, only two points $(x,y)=(0,\haf\pi)$ and
$(\haf\pi,0)$ are the stationary points of the phase and hence
only these contribute to the leading term of \eqref{invar-bound}.
Moreover this contribution was computed in the course of the proof
of Theorem \ref{invar-thm}. This contribution is coming from an
integral of $v$ over the neighborhood of the size smaller than
$|\lm_i|^{-\haf+\eps}$ for any $\eps>0$. This again follows from
the stationary phase method. The function $v$ is well approximated
on a small enough interval by its value in the center of this
interval. Hence by letting $r\to\8$ but keeping it smaller than
$|\lm_i|^{\haf-\eps'}$ with $\eps'>\eps$ we see that for any
smooth function $\psi$ the value of the integral against $\rho_r$
has the same leading term as the integral against $dU_i$ as
$|\lm_i|\to\8$ and $r\to\8$. The constant $C_{\psi,\eps}$ in the
$O$-term is bounded by appropriate derivative of $v$ at stationary
points.

\end{proof}

\section{Spectral localization of eigenfunctions under the action of PDO}
\subsection{Spectral localization} We consider now more general
matrix coefficients. Let $a\in C^\8(X)$ be a symbol and $Op(a)$
the corresponding pseudo-differential operator. We assume for
simplicity that $a$ belongs to an automorphic representation of
class one. For a fixed symbol $a$ we are interested in the
decomposition of $Op(a)\phi_i$ with respect to the basis of
eigenfunctions $\{\phi_j\}$ as $|\lm_i|\to\8$.

\begin{thm}{prop} For for a fixed symbol $a\in C^\8(X)$ and for
any $N>0$ the following bound holds
\begin{align}
\label{phi-i-j-N} |<Op(a)\phi_i,\phi_j>|= O_N(|\lm_i-\lm_j|^{-N})\
\end{align}
with the constant in the $O$-term depends on $N$ and on the
symbol.
\end{thm}
\begin{proof}
To prove \eqref{phi-i-j-N} we need to analyze the values of
\begin{align}
l_{\pi_\mu,\pi_i,\pi_j}^{aut}(a,e_{\pi_i},\phi_j)=
a_{\pi_\mu,\pi_{\lm_i},\pi_{\lm_j}}
l_{\pi_\mu,\pi_i,\pi_j}^{mod}(a,e_{\pi_i},\phi_j)\ .
\end{align}
We saw that coefficients $a_{\pi_\mu,\pi_{\lm_i},\pi_{\lm_j}}$ are
polynomially bounded in $\lm_i$. On the other hand it follows from
the stationary phase method that the structure of the model
trilinear invariant functionals $l^{mod}$ is governed by the
presence of critical points of the phase and singularities of the
amplitude. It is easy to see that as  $|\lm_i-\lm_j|\to\8$ the
phase function in the kernel of $l^{mod}$ does not have critical
points with respect to $x$-integration in \eqref{integ-reduc1} and
non-degenerate critical point with respect to $z$-integration. The
amplitude of the kernel is becoming a smooth function after
integration against the smooth function $a$. Hence from the
stationary phase method and the bound \eqref{coef-a-bound} on
coefficients $a_{\pi_\mu,\pi_{\lm_i},\pi_{\lm_j}}$ we obtain the
bound claimed.
\end{proof}
\subsection{Conjectural density}
According to the proposition above we see that for a
fixed symbol $a$ the spectral density of $Op(a)\phi_i$ is
essentially supported on a very short interval around $\lm_i$
itself. Hence, under the action of $Op(a)$ the eigenfunctions
$\phi_i$ are spectrally localized in short intervals.

This makes a question about spectral density of $Op(a)\phi_i$
inside the interval $|\lm_i-\lm_j|\ll|\lm_i|^\eps$ interesting. We
note that this question was also raised in a connection with the
quantum unique ergodicity conjecture. We may speculate about the
size of coefficients $<Op(a)\phi_i,\phi_j>$ on the basis of a
conjecture similar to Conjecture \ref{CONJ}. One is lead to
conjecture (though, solely on the basis of examples of arithmetic
surfaces, see \cite{Sa2}) that the coefficients
$a_{\pi_\mu,\pi_{\lm_i},\pi_{\lm_j}}$ are of the order
$\max\{|\lm_i|^\eps,|\lm_j|^\eps\}$ for any $\eps>0$. On the other
hand it is also expected that these coefficients are not small on
the average (though some of them could be zero). Namely, one
expects, for example, that for all $\lm_i$ and any fixed $B>0$ the
lower bound
\begin{align}\label{a-average}
\sum_{|\lm_i-\lm_j|\leq B}
|a_{\pi_\mu,\pi_{\lm_i},\pi_{\lm_i}}|^2\geq c|\lm_i| \ ,
\end{align}
holds for some $c>0$. This again is consistent with the
Lindel\"{o}f conjecture since according to the  Weyl law the
number of terms in the sum above is of order $|\lm_i|$.

On the other hand it is easy to see from the stationary phase
method that
\begin{align}\label{l-bound}
|l_{\pi_\mu,\pi_i,\pi_j}^{mod}(a,e_{\pi_i},\phi_j)|\asymp|\lm_i|^{-\haf}
\end{align}
for $|\lm_i-\lm_j|\leq B$. We deduce from this that for a fixed
symbol $a$ the spectral density of $Op(a)\phi_i$ is supported in
the interval $|\lm_i-\lm|\ll|\lm_i|^\eps$ (as we have shown in
Theorem \ref{prop}) and conjecture that it has the absolute value
of order $|\lm_i|^{-\haf+\eps}$ at most on this interval. Hence,
$Op(a)$ spreads $\phi_i$ evenly on this interval.

Similarly, the conjectural upper bound for the coefficients
$a_{\pi_\mu,\pi_{\lm_i}\pi_{\lm_j}}$ implies that the matrix
coefficients satisfy
\begin{align}
\label{phi-i-j}
|<Op(a)\phi_i,\phi_j>|\ll |\lm_i|^{-\haf}\ ,
\end{align}
as $|\lm_i-\lm_j|\to 0$ and $|\lm_i|\to\8$.


\end{document}